\documentclass[letterpaper, 10pt, journal, twoside]{IEEEtran}      

\IEEEoverridecommandlockouts                             

\usepackage{graphics} 
\usepackage{epsfig} 
\usepackage{times} 
\usepackage{amsmath} 
\usepackage{amssymb}  
\usepackage{tikz}
\usepackage{mathtools}
\usepackage{comment}
\usepackage{tikz}
\usepackage{pgfplots} 
\usepackage{algorithm}
\usepackage{algorithmic}

\allowdisplaybreaks[0]

\title{Fast Optimal Energy Management with Engine On/Off Decisions for Plug-in Hybrid Electric Vehicles
}

\author{Sebastian East and Mark Cannon
\thanks{S.~East is with the Department of Engineering Science, University of Oxford, OX1 3JP, {\tt\footnotesize sebastian.east@eng.ox.ac.uk}, (+44) 1865 273000.}%
\thanks{M.~Cannon is with the Department of Engineering Science, University of Oxford, OX1 3JP, {\tt\footnotesize mark.cannon@eng.ox.ac.uk}, (+44) 1865 273000.}%
}

\begin{document}

\maketitle
\thispagestyle{empty}
\pagestyle{empty}

\begin{abstract}

In this paper we demonstrate a novel alternating direction method of multipliers (ADMM) algorithm for the solution of the hybrid vehicle energy management problem considering both power split and engine on/off decisions. The solution of a convex relaxation of the problem is used to initialize the optimization, which is necessarily nonconvex, and whilst only local convergence can be guaranteed, it is demonstrated that the algorithm will terminate with the optimal power split for the given engine switching sequence. The algorithm is compared in simulation against a charge-depleting/charge-sustaining (CDCS) strategy and dynamic programming (DP) using real world driver behaviour data, and it is demonstrated that the algorithm achieves 90\% of the fuel savings obtained using DP with a 3000-fold reduction in computational time. 

\end{abstract}


\section{INTRODUCTION}
\IEEEPARstart{O}{ne} of the current practical limitations of all-electric vehicles is the restricted driving range available from a single charge due to the low-energy density of battery technology. Plug-in hybrid electric vehicles (PHEVs), where a modestly sized electrical propulsion system is complemented with an internal combustion engine, are a common compromise \cite{EhsaniTextbook}. A problem associated with hybrid powertrains, however, is that the additional power source requires that the fraction of power delivered from the electrical and internal combustion systems is actively controlled by the vehicle. It has been demonstrated that by modulating this fraction throughout the journey the total fuel consumption can be reduced, and this is known as the `energy management problem' \cite{Sciarretta2007}.

If the future driver behaviour can be determined with a high degree of certainty it is possible to formulate the energy management problem as an optimal control problem, and a wide range of optimization based methods have been demonstrated for its solution (for a recent survey see \cite{MarinaMartinez2016}). Dynamic programming (DP) can guarantee global optimality when the powertrain dynamics are considered without simplification, but the computational complexity of DP increases exponentially with the number of decision variables \cite{Sundstr2009} and is therefore not feasible for an online implementation or for use as part of a model predictive control (MPC) framework.

Convex optimization approaches to the problem have been investigated to reduce the computational cost, and it has been demonstrated that if the engine switching and gear selection are determined by a heuristic or prior optimization routine, the optimal power split can be determined from the solution of a convex optimization problem \cite{Egardt2014}. Introducing engine switching and gear selection decisions into the optimization problem is a significant challenge, however, as both require integer decision variables. A range of approaches have been taken, including mixed integer programming \cite{Reeven2012}, Pontryagin's Minimum Principle (PMP) \cite{Schmid2018}, convex relaxation \cite{Josevski2016a}, and genetic algorithms \cite{Guo2018}. In \cite{Elbert2014} it was shown in simulations that the globally optimal solution can be obtained by alternating between a convex problem to determine the optimal power split for a fixed engine switching sequence, and PMP for the optimal engine switching given a fixed power split, and that this is generally faster than dynamic programming. However we note that convergence of this approach to the globally optimal solution cannot be proven in general. A similar approach was taken in \cite{Nuesch2014} to also include the gear selection and engine switching cost.

Whilst these approaches have demonstrated improved performance compared to DP with minimal model reduction, multiple iterations of the convex power split sub-problem still incur a significant computational burden, and they are too slow for an online solution: up to 155s were required in \cite{Nuesch2014}. Here, we instead use an alternating direction method of multipliers (ADMM) \cite{Wang2014} algorithm to first obtain the optimal control values for a convex relaxation of the problem, and then use these values to initialise an ADMM solution of the nonconvex problem. Whilst convergence can only be guaranteed to a point that satisfies first order (neccessary) optimality conditions, the algorithm's structure ensures that the power split is optimal for the engine switching control at termination. The performance of the algorithm is compared through simulation on real driver data w.r.t. a charge-sustaining/charge-depleting (CDCS) strategy and DP.

This paper extends our earlier work on the PHEV energy management problem. In \cite{Buerger} we demonstrated a dual-loop MPC framework with a projected-Newton method, but only considered a terminal state of charge constraint and ignored engine switching. In \cite{EastCDC2018} we proved that a formulation including general state of charge limits is convex when stated as a function of the battery power, and proposed an ADMM algorithm for its solution, but still did not consider engine switching. In \cite{East2018a} we proved that the algorithm of \cite{EastCDC2018} is convergent even in the presence of nonconvexity, but only to a point that satisfies the first order (neccessary) conditions. 

ADMM was previously used in \cite{Takapoui2017} as a heuristic for energy management with engine switching, but all losses in the electric system were ignored and the algorithm was not necessarily convergent. Furthermore, all hard limits on power and state of charge were ignored, so it is not suitable for embedded control. 

The rest of the paper is structured as follows: in section 2 the energy management problem is mathematically formalised and in section 3 an ADMM algorithm is presented for its solution. Section 4 details the numerical studies, and the paper is concluded in section 5.

\section{Problem Formulation}\label{section_2}

For the work presented here, we extend the model previously described in \cite{EastCDC2018} to include engine switching, and in places we refer the reader to that paper for more details. 

Figure 1 shows a simplified diagram of a parallel PHEV powertrain and is used to define the power transfers, $P$. It is assumed that the clutch engagement and engine state are coupled, and both given at the $k$th sampling instant by $\sigma_k$, where $\sigma_k = 1$ indicates that the engine is on and the clutch engaged, and $\sigma_k = 0$ indicates that the engine is off and the clutch disengaged.
\begin{figure}
\centering
\scalebox{0.72}{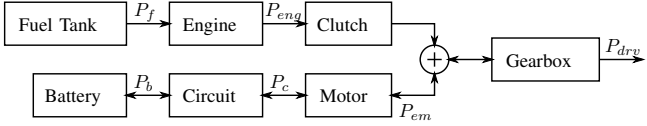}
\caption{Powertrain Model (see \protect\cite{EastCDC2018} for full interpretation of each variable) }
\label{fig_powertrain_model}
\end{figure} 
It is assumed that an accurate prediction of future velocity, $v$, and road gradient, $\theta$, is available over a prediction horizon of $N$ steps, from which the power demand, $P_{drv}$, is given for $k = 0,\dots,N-1$ by
\begin{multline*}
P_{drv,k} = ( m \dot{v}_k + \frac{1}{2} \rho_a v_k^2 C_d A \\
+ C_r mg \cos \theta_k + m g \sin \theta_k ) v_k - P_{brk,k},
\end{multline*}
where $\dot{v}_k$ can be obtained from numerical differentiation, $m$ is the vehicle mass, $\rho_a$ is the density of air, $C_d$ is the drag coefficient, $A$ is the frontal area, $C_r$ is the rolling resistance, $g$ is the acceleration due to gravity, and $P_{brk,k}$ is the mechanical braking power. It is assumed that the mechanical brake input is determined by an external controller so the total demand power is defined at all times. The angular speed of the drivetrain input, $\omega_d$, is given at all times by a simple gear switching heuristic function based on vehicle velocity e.g. $\omega_{d,k} = f(v_k)$ (see \cite{EastCDC2018}). It is assumed that when the drivetrain speed is below the minimum engine speed the clutch is disengaged and the engine is off, and that a negative demand power is delivered from the motor only. These assumptions are used to define the sets
\begin{align*}
\mathcal{P} &= \{k: P_{drv,k} \geq 0, \omega_{d,k} \geq \underline{\omega}_{eng} \}, \\
\mathcal{B} &= \{k : P_{drv,k} < 0, \omega_{d,k} \geq \underline{\omega}_{eng} \}, \\
\mathcal{C} &= \{k: \omega_{d,k} < \underline{\omega}_{eng} \},
\end{align*}
so that power split and engine state need to be determined for $k \in \mathcal{P}$, and only engine swtiching needs to be determined for $k \in \mathcal{B}$, i.e
\begin{align*}
P_{drv,k} &= \begin{cases} P_{eng,k} + P_{em,k} & k \in \mathcal{P} \\ P_{em,k} & k \notin \mathcal{P} \end{cases}, \\
\sigma_k &\in \begin{cases} \{0,1\} & k \notin \mathcal{C} \\ 0 & k \in \mathcal{C} \end{cases}.
\end{align*}

 The rotational speed of the engine, $\omega_{eng}$, and motor, $\omega_{em}$, are therefore given by $\omega_{em,k} = \omega_{d,k}$ and $\omega_{eng,k} = \sigma_k \omega_{d,k}$. Starting from the assumption that the engine is on and engaged for $k\notin \mathcal{C}$, engine loss functions, $f_k$, and motor loss functions, $h_k$, can be obtained $\forall \ k$ in the convex quadratic form
\begin{align*}
f_k(P_{eng,k}) &= \alpha_{2,k} P_{eng,k}^2 + \alpha_{1,k} P_{eng,k} + \alpha_{0,k}, \\
h_k(P_{em,k}) & = \beta_{2,k} P_{em,k}^2 + \beta_{1,k} P_{em,k} + \beta_{0,k},
\end{align*}
where $\alpha_{2,k},\beta_{2,k} \geq 0$. Losses in the battery are modelled using the equivalent circuit model
$$
g_k(P_{em,k}) = \frac{V^2}{2R} \left( 1 - \sqrt{1 - \frac{4 R}{V^2} h_k (P_{em,k})} \right),
$$
where it is assumed that the open circuit voltage, $V$, and resistance, $R$, are constant. 

Limits are introduced on both motor and engine power for $k \in \mathcal{P}$ to ensure that their respective torque limits are not exceeded, and to ensure that $g_k$ and $f_k$ are all non-decreasing (this is reasonable as it would be expected that an increase in output power would require an increase in input power) and real-valued (see \cite{EastCDC2018}). These can then be consolidated into upper and lower limits on motor power, $\overline{P}_{em,k}$ and $\underline{P}_{em,k}$, so that limits on battery power can be given $\forall \ k$ as
\begin{align*}
\underline{P}_{b,k} = \begin{cases} g_k(\underline{P}_{em,k} ) & k \in \mathcal{P} \\ g_k(P_{drv,k}) & k \notin \mathcal{P} \end{cases} \ \overline{P}_{b,k} = \begin{cases} g_k(\overline{P}_{em,k} ) & k \in \mathcal{P} \\ g_k(P_{drv,k}) & k \notin \mathcal{P} \end{cases}.
\end{align*}
The restriction of $g_k$ and $f_k$ allows the inverse function of $g_k$ to be defined as
$$
g_k^{-1}(P_{b,k}) = -\frac{\beta_{1,k}}{2 \beta_{2,k}} + \sqrt{-\frac{R P_{b,k}^2}{\beta_{2,k} V^2} + \frac{P_{b,k} - \beta_{0,k}}{\beta_{2,k}} + \frac{\beta_{1,k}^2}{4 \beta_{2,k}^2}}.
$$
The dynamics of clutch engagement and engine switching are ignored, and it is assumed that the engine switching variable $\sigma_k$ introduces two effects: A) all of the power is delivered by the electric system when the engine is off B) the fuel consumption when the engine is off is zero. Effect A) is introduced with the constraint
\begin{equation}\label{eqn_constraint}
g_k(P_{drv,k}) + \sigma_k \gamma_k \leq P_{b,k} \leq g_k(P_{drv,k}) + \sigma_k \delta_k,
\end{equation}
where $\gamma_k = \underline{P}_{b,k} - g_k(P_{drv,k})$ and $\delta_k = \overline{P}_{b,k} - g_k(P_{drv,k})$.
This constraint is convex subject to the condition that $\gamma_k \leq \delta_k$, which is true whenever the problem is feasible, i.e $\overline{P}_{b,k} \geq g_k(P_{drv,k}) \geq \overline{P}_{b,k}$. Effect B) can be introduced by defining the fuel consumption with the function
\begin{align*}
\dot{m}_f(P_b,\sigma) =& \sum_{k \in \mathcal{P}} \Big[ f_k(P_{drv,k} - g_k^{-1}(P_{b,k})) \\
&+ (\sigma_k - 1)f_k (0) \Big] + \sum_{k \in \mathcal{B}} \sigma_k f_k(0),
\end{align*}
as we note from (\ref{eqn_constraint}) that $\sigma_k = 0 \Rightarrow f_k(P_{drv,k} - g_k^{-1}(P_{b,k})) = f_k(0)$. It is demonstrated in \cite{EastCDC2018} that $f_k(P_{drv,k} - g_k^{-1}(P_{b,k}))$ is convex and non-increasing, so the function $\dot{m}_f(P_b,\sigma)$ is therefore convex and separable w.r.t both $\{ \sigma,P_b\}$ and the individual elements $\sigma_k$ and $P_{b,k}$. Finally, a `driveability' cost is introduced as
$$
d(\sigma) = \frac{k_d}{2} \sigma^\top(\Psi^{-1})^\top \Psi^{-1} \sigma,
$$
where $\Psi$ is a lower triangular matrix of ones of the conformal size, and the value of $k_d$ can either be chosen to accurately model the fuel consumed when switching the engine on, or chosen arbitrarily to prevent the engine control from chattering. 
The state of charge dynamics are approximated using Euler integration (with a  sampling interval of 1s) as $E_{k+1} = E_k - P_{b,k}$, and the state of charge of the battery is subject to constant upper and lower bounds, $\overline{E}$ and $\underline{E}$.
The optimization problem to obtain the control inputs, $P_{b}^\star$ and $\sigma^\star$, that minimize fuel consumption is therefore given by
\begin{equation}\label{eqn_optimization_problem}
\begin{aligned}
\min_{P_b,\sigma} \ & \sum_{k \in \mathcal{P}} \Big[ f_k(P_{drv,k} - g_k^{-1}(P_{b,k})) + (\sigma - 1)f_k (0) \Big]\\
 &+ \sum_{k \in \mathcal{B}} \sigma_k f_k(0) + \frac{k_d}{2} \sigma^\top(\Psi^{-1})^\top \Psi^{-1} \sigma \\
\text{s.t.} \ &E = \Phi E_0 - \Psi P_b \\
& \begin{rcases} \underline{E} \leq E_k \leq \overline{E} \\
g_k(P_{drv,k}) + \sigma_k \gamma_k \leq P_{b,k} \leq g_k(P_{drv,k}) + \sigma_k \delta_k \end{rcases} \forall k \\
&\sigma_k \in \begin{cases} \{0 \} & k \in \mathcal{C} \\ \mathcal{S} & k \notin \mathcal{C} \end{cases}.
\end{aligned}
\end{equation}
From the properties of the cost functions and constraints demonstrated above, this problem is nonconvex when $\mathcal{S} = \{0,1\}$, but becomes convex if the set is relaxed to $\mathcal{S} = [0,1]$. If $\sigma_k = 0 \ \forall \ k$ and $P_{b,k} = g_k(P_{drv,k}) \ \forall \ k$ (i.e, the engine is off for the entire journey) satisfy the state constraints then clearly these values are optimal. Otherwise, an algorithm is required for its solution.
\section{Alternating Direction Method of Multipliers}

\begin{figure}
\begin{center}
\begin{tikzpicture}
\begin{axis}[scale = 0.8, axis lines = center, enlargelimits = true, clip = false, xlabel = {$\sigma_k$}, ylabel = {$\eta_k$}, xmajorticks = false, xticklabels = {,,}, ytick = {0.5}, yticklabels = {$g_k(P_{drv,k})$}]
\addplot [black, mark = none] coordinates {(-0.5,0) (1.5,2)} node [right] {$g_k(P_{drv,k}) + \sigma_k \delta_k$};
\addplot [black] coordinates {(-0.5,1) (1.5,-1)} node [right] {$g_k(P_{drv,k}) + \sigma_k \gamma_k$};
\addplot [black] coordinates {(1,-1) (1,2)} node [above] {$\sigma_k = 1$};
\draw [fill=blue,opacity=0.1] (axis cs:0,0.5) -- (axis cs:1,-0.5) -- (axis cs:1,1.5);
\addplot[only marks, red] coordinates {(0.66,0.5)} node [right] {1};
\addplot[only marks, red] coordinates {(0,0.5)} node [above right] {2};
\addplot[only marks, red] coordinates {(0.5,1)} node [above left] {6};
\addplot[only marks, red] coordinates {(1,1.5)} node [above left] {3};
\addplot[only marks, red] coordinates {(1,0.5)} node [right] {5};
\addplot[only marks, red] coordinates {(1,-0.5)} node [below left] {4};
\addplot[only marks, red] coordinates {(0.5,0)} node [below left] {7};
\end{axis}
\end{tikzpicture}
\caption{Constraint set defined by $h^{\mathcal{S}}_k(\sigma_k) + h^{P_b}_k(\eta_k,\sigma_k)$ for $\mathcal{S} = [0,1]$. \label{fig_etaupdate}}
\end{center}
\end{figure}
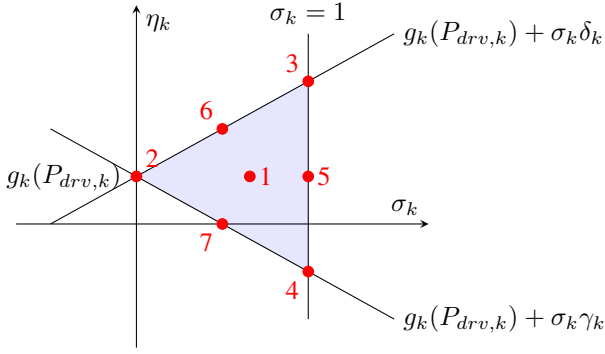

We propose an alternating direction method of multipliers algorithm for the solution of (\ref{eqn_optimization_problem}) in two phases: firstly, it is used to find the global solution to the convex relaxation of the problem with $\mathcal{S} = [0,1]$, then this solution is used as the initialisation for the problem with $\mathcal{S} = \{0,1\}$. Two dummy variables, $\eta$ and $\zeta$, are introduced (their inclusion is explained below), and (\ref{eqn_optimization_problem}) is then given as the equality constrained problem
\begin{equation}\label{eqn_optimization_problem2}
\begin{aligned}
\min_{P_b,\sigma} \ & \dot{m}_f(P_b,\sigma) + d(\kappa) + H^\mathcal{S}(\sigma) + H^E(E) + H^{P_b}(\eta, \sigma) \\
\text{s.t.} \ &E = \Phi E_0 - \Psi \zeta \\
& \zeta =  P_b = \eta \\
& \kappa = \sigma
\end{aligned}
\end{equation}
with the indicator functions
\begin{align*}
H^\mathcal{S}(\sigma) =& \sum_{k=0}^{N-1} h^\mathcal{S}(\sigma_k), \ h^\mathcal{S}(\sigma_k) = \begin{cases} 0 & \sigma_k \in \mathcal{S} \\ \infty & \text{otherwise} \end{cases} \\
H^E(E) =& \sum_{k=1}^{N} h^E(E_{k}), \ h^E(E) = \begin{cases} 0 & \underline{E} \leq E_{k} \leq \overline{E} \\ \infty & \text{otherwise} \end{cases} \\
H^{P_b}(
\eta,\sigma) =& \sum_{k=0}^{N-1} h^{P_b}_k(\eta_k,\sigma_k), \\
& h^{P_b}_k(\eta_k,\sigma_k) = \begin{cases} \infty & g_k(P_{drv,k}) + \sigma_k \gamma_k \geq \eta_k \\
\infty & \eta_k \geq g_k(P_{drv,k}) + \sigma_k \delta_k  \\ 0 & \text{else} \end{cases}.
\end{align*}
The augmented Lagrangian associated with (\ref{eqn_optimization_problem2}) is given by
\begin{equation}\label{eqn_lagrangian}
\begin{aligned}
& \mathcal{L} (\kappa, P_b, E, \sigma, \eta, \zeta, \lambda_1, \lambda_2, \lambda_3, \lambda_4) \\
= & \ \dot{m}_f(P_b,\sigma) + d(\kappa) + H^{\mathcal{S}}(\sigma) + H^E(E) + H^{P_b}(\eta,\sigma) \\
& + \frac{\rho_1}{2} \| \Phi E_0 - \Psi \zeta - E + \lambda_1 \|^2  + \frac{\rho_2}{2} \|P_b - \zeta + \lambda_2 \|^2  \\ 
&+  \frac{\rho_3}{2} \|P_b - \eta + \lambda_3 \|^2 + \frac{\rho_4}{2} \| \kappa - \sigma + \lambda_4 \|^2,
\end{aligned}
\end{equation}
and after initialising the variables with 
\begin{align*}
\sigma^0 &= 0, \quad \eta^0_k = \zeta^0_k = P_{b,k}^0 = g_k(P_{drv,k}), \\ 
E^0 &= \Pi_{[\underline{E},\overline{E}]}(\Phi E_0 - \Psi \zeta^0), \quad \lambda_1^0 = \lambda_2^0 = \lambda_3^0 = \lambda_4^0 = 0,
\end{align*}
where
\begin{align*}
\Pi_{[\underline{x},\overline{x}]}(x) &= \{\pi_{[\underline{x}_1,\overline{x}_1]}(x_1),\dots,\pi_{[\underline{x}_N,\overline{x}_N]}(x_N) \}, \\
\pi_{[\underline{x}_i,\overline{x}_i]}(x_i) &= \min \{ \overline{x}_i, \max \{ \underline{x}_i, x_i \} \},
\end{align*}
the ADMM iteration is obtained by minimising the augmented Lagrangian w.r.t each of the variables in turn, in the order of the arguments of $\mathcal{L}$ in (\ref{eqn_lagrangian}). Firstly, at iteration $j$, the $\kappa$ update has the analytical solution
\begin{align*}
\kappa^{j+1} =& \big(k_d(\Psi^{-1})^\top \Psi^{-1} + \rho_4 I \big)^{-1} \rho_4 (\sigma^j - \lambda_4^j),
\end{align*}
then the $P_{b,k}$ update for $k \in \mathcal{P}$ is found from 
\begin{multline*}
P^{j+1}_{b,k} = \arg \min_{P_{b,k}} ( f_k ( P_{drv,k} - g_k^{-1}(P_{b,k})) \\
+ \frac{\rho_2}{2} (P_{b,k} - \zeta_k^j + \lambda_{2,k}^j )^2 + \frac{\rho_3}{2}(P_{b,k} -\eta^{j}_k + \lambda_{3,k}^j)^2 ).
\end{multline*}
This is an unconstrained convex optimization problem that is solved using a Newton method with a backtracking line search (as in \cite{EastCDC2018}). For $k \notin \mathcal{P}$ this minimization has the analytical solution
$$
P_{b,k}^{j+1} = \frac{\rho_2 ( \zeta_k^j - \lambda_{2,k}^j) + \rho_3(\eta_k^j - \lambda_{3,k}^j)}{\rho_2 + \rho_3}.
$$
The $E$ update is 
$$
E^{j+1} = \Pi_{[\underline{E}, \overline{E}]} (\Phi E_0 - \Psi \zeta^{j} + \lambda_1^j).
$$
The variables $\eta_k$ and $\sigma_k$ are simultaneously updated for $k \in \mathcal{P} \cup \mathcal{B}$ as
\begin{multline}\label{eqn_etasigmaupdate}
( \eta_k,\sigma_k )^{j+1} = \arg \min_{\eta_k,\sigma_k} \Big( \sigma_k f_k(0) + \frac{\rho_3}{2}(P_{b,k}^{j+1} - \eta_k + \lambda_{3,k}^j)^2 \\
+ \frac{\rho_4}{2}(\kappa_k^{j+1} - \sigma_k + \lambda_{4,k}^j)^2 + h^{\mathcal{S}}_k(\sigma_k) + h^{P_b}_k(\eta_k,\sigma_k) \Big)
\end{multline}
(note that whilst the first term should be $(\sigma_k - 1)f_k(0)$ for $k \in \mathcal{P}$, the inclusion of a constant does not affect the solution). When $\mathcal{S} = [0,1]$ this is a convex inequality-constrained quadratic program for which there are 7 possible solutions as shown in Figure \ref{fig_etaupdate}: the unconstrained minimum ({\color{red}1}), three at the vertices ({\color{red}2-4}), and three on the edges ({\color{red}5-7}). Each of these has an analytical solution, so we obtain the update by finding each possible solution and determining which minimises (\ref{eqn_etasigmaupdate}). When $\mathcal{S} = \{0, 1\}$ the constraints become nonconvex, but the solution is still simple to obtain: either $(\eta_k,\sigma_k)^{j+1} = (g(P_{drv,k}),0)$, or $(\eta_k,\sigma_k)^{j+1} = (\pi_{[\underline{P}_{b,k},\overline{P}_{b,k}]}(P_{b,k}^{j+1} + \lambda_{3,k}^j),1)$, so we find both and determine which is the minimizing argument of (\ref{eqn_etasigmaupdate}). For $k \in \mathcal{C}$ the update is simply given by $(\eta_k,\sigma_k)^{j+1} = (g(P_{drv,k}),0)$.
Finally, the remaining updates have the analytical solutions
\begin{alignat*}{2}
\zeta^{j+1} =& \big( \rho_2 I + \rho_1 \Psi^\top \Psi \big)^{-1} \\
& \big( \rho_2 (P_{b}^{j+1} + \lambda_2^j) + \rho_1 \Psi^\top (\Phi E_0 - E^{j+1} + \lambda_1^j) \big), \\
\lambda_1^{j+1} =& \lambda_1^j + \Phi E_0 - \Psi \zeta^{j+1} - E^{j+1}, \\
\lambda_2^{j+1} =& \lambda_2^j + P_b^{j+1} - \zeta^{j+1}, \\
\lambda_3^{j+1} =& \lambda_3^j + P_b^{j+1} - \eta^{j+1}, \\
\lambda_4^{j+1} =& \lambda_4^j + \kappa^{j+1} - \sigma^{j+1}.
\end{alignat*}
This is not the only possible ADMM formulation: additional dummy variables are included to reduce the complexity of each variable update, but increase the total number of variable updates per iteration. We include $\kappa$ so that the minimisation of the switching cost, $d(\sigma)$, can be performed analytically, and $\zeta$ so that the constraints on $E$ can be applied element-wise. $\eta$ is included so that the combined $\eta,\sigma$ update has quadratic cost and can be performed cheaply using the rules demonstrated in Figure \ref{fig_etaupdate}; if the function $f_k(P_{drv,k} - g^{-1}_k(P_{b,k}))$ is approximated as a quadratic function, $\eta$ is therefore not required. 

\subsection{Convergence, Optimality, \& Complexity}

We first rearrange problem (\ref{eqn_optimization_problem2}) into the form
\begin{align*}
\min_{u,x} &\quad \hat{f}(u) + \hat{g}(x) \\
\text{s.t} \ & \quad Au + Bx = c,
\end{align*}
where
\begin{align*}
u =& (\kappa, P_b, E), \quad x = (\sigma, \eta, \zeta),  \\
\hat{f}(u) =& \sum_{k \in \mathcal{P}} f_k(P_{drv} - g^{-1}(P_{b,k})) + d(\kappa)  + H^E(E), \\
\hat{g}(x) =& \sum_{k \in \mathcal{P} \cup \mathcal{B}} \sigma_k f_k(0) + H^\mathcal{S}(\sigma) + H^{P_b}(\eta, \sigma),  \\
A =& \begin{bsmallmatrix}
0 & 0 & -I \\ 0 & I & 0 \\ 0 & I & 0 \\ I & 0 & 0 
\end{bsmallmatrix}, \ B = \begin{bsmallmatrix}
0 & 0 & -\Psi \\ 0 & 0 & -I \\ 0 & -I & 0 \\ -I & 0 & 0
\end{bsmallmatrix} \ c = (-\Phi E_0,0,0,0). 
\end{align*}
This is the form used in \cite{Wang2014}, where it is demonstrated that, for a convex problem, an equivalent ADMM iteration will converge to the optimal solution as the residuals defined by $r^{j+1} = Au^{j+1} + Bx^{j+1} - c$ and $s^{j+1} = A^\top \rho B (x^{j+1} - x^j)$ converge to zero, where $\rho = \text{diag}\{\rho_1 I, \rho_2I, \rho_3I, \rho_4I \} \succ 0 $. For the convex initialisation phase of the algorithm in which $\mathcal{S} = [0,1]$, the values of the decision variables at termination can therefore be made arbitrarily close to the globally optimal values by setting the value of $\epsilon$ arbitrarily close to zero, and using the termination criterion $\max \left\{ \|r^{j+1}\|, \|s^{j+1} \| \right\} \leq \epsilon$. 

We have previously proved in \cite{East2018a} that the residuals will still converge to zero if the problem is nonconvex, and the algorithm will therefore converge to a point satisfying the first order optimality conditions when $\mathcal{S} = \{0,1\}$ (the proof is for non-affine equality constraints, but an identical argument holds for integer decision variables). The principle of the approach presented here is that the solution of the convex relaxation is used to initialise the algorithm in the nonconvex case. Furthermore, we note that the nonconvexity is only associated with $\sigma$, and that if the value of $\sigma$ were fixed, problem (\ref{eqn_optimization_problem2}) would become convex. This implies that for any sequence of iterates in which $\sigma$ is stationary, the other variables will have converged to their corresponding optimal values, which therefore also implies that the power split at termination will be optimal for the engine switching sequence obtained (this result is as strong as that shown in \cite{Elbert2014}). The algorithm is summarised in Algorithm \ref{algorithm_ADMM}.
\begin{algorithm}
\caption{ADMM\label{algorithm_ADMM}}
\begin{algorithmic}[1]
\STATE Initialise $\{\kappa,P_b,\sigma,\eta,\zeta,E,\lambda_{1-4}\}^0$, $j \gets 0$, $\mathcal{S} \gets [0,1]$
\WHILE {$\max \left\{ \|r^{j+1}\|, \|s^{j+1} \| \right\} > \epsilon$}
\STATE Calculate $\{\kappa,P_b,\sigma,\eta,\zeta,E,\lambda_{1-4}\}^{j+1}$
\STATE $j \gets j+1$
\ENDWHILE
\STATE $\mathcal{S} \gets \{0,1\}$
\WHILE {$\max \left\{ \|r^{j+1}\|, \|s^{j+1} \| \right\} > \epsilon$}
\STATE Calculate $\{\kappa,P_b,\sigma,\eta,\zeta,E,\lambda_{1-4}\}^{j+1}$
\STATE $j \gets j+1$
\ENDWHILE
\end{algorithmic}
\end{algorithm}

The computation of each iteration is dominated by the $\kappa$, $P_b$, $\sigma$, $\eta$, and $\zeta$ updates. All other updates (including $r$ and $s$ when analytically block-multiplied) scale linearly with $N$ (multiplication by $\Psi$ is equivalent to a cumulative sum). The $P_b$, $\sigma$ and $\eta$ updates also scale linearly with $N$ when implemented sequentially, or are constant if each $k$ element is updated in parallel and sufficient threads are available. The dense matrix inversions associated with the $\kappa$ and $\zeta$ updates can be computed offline as they do not include any decision variables, so only dense matrix-vector multiplications are required, and the computational complexity of each iteration is therefore $\mathcal{O}(N^2)$. Additionally, the updates for $\kappa$ and $\zeta$ consist of multiplications by Toeplitz matrices that can be implemented as (stable) linear filtering operations with a storage requirement of $\mathcal{O}(N)$.

\section{Numerical Studies}

The performance of the ADMM algorithm was investigated through simulation in comparison with both a CDCS strategy and an approximately optimal DP implementation (the algorithm of \cite{Takapoui2017} was not included as the lack of hard limits on state of charge mean that it cannot be compared in any meaningful way). For both optimization-based strategies it was assumed that the future driver behaviour was known with complete precision. A simple CDCS strategy was assumed where the engine was switched off and all power was delivered from the motor until the lower state constraint was violated, after which the engine was permanently switched on, and used to provide all of the positive demand power whenever the state of charge was below its lower constraint. The DP algorithm was modified from that presented in \cite{EastCDC2018} to include the engine switching control variable and engine switching cost, with the state of charge discretised to 0.1\% intervals and the battery power control input discretised to 1\% intervals. The values $\rho_1 = 8.86\times10^{-9}$ and $\rho_2 = 2.34\times10^{-4}$ were taken from \cite{EastCDC2018}, and we set $\rho_3 = \rho_2$. Also, $\rho_4$ was set at 2$\times10^3$ after using a parameter sweep similar to that detailed in \cite{EastCDC2018}. The termination threshold, $\epsilon$, was set at $7\times10^4$ for both the initial convex phase and the subsequent nonconvex phase of the ADMM algorithm, and $k_d$ was set arbitrarily at $10^4$.

\begin{figure}
\begin{center}
\includegraphics[scale=1]{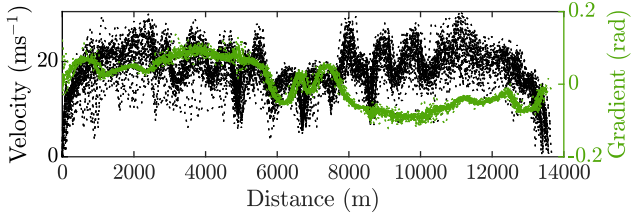}
\caption{Velocity and gradient data against distance. \label{fig_velocity}}
\end{center}
\end{figure}

The velocity and road gradient data used to generate the power demand profiles is shown in Figure \ref{fig_velocity}. This is real drive-test data taken from 49 instances of a single $\sim$13km route driven by four different drivers. The method detailed in section \ref{section_2} was used to obtain the demand power for this data, where it was assumed that the mechanical brake provided none of the braking power. The vehicle was modelled as a 1800kg passenger vehicle with a 100kW petrol internal combustion engine, a 50kW electric motor, and a 21.5Ah lithium-ion battery with a 350V and 0.1$\Omega$ open circuit voltage and resistance. The battery was initialised at 60\% and constrained to between 40\% and 70\% to ensure that the state constraints were strongly active at the solution. The simulations were programmed in Matlab on a 2.6GHz Intel Core i7-6700HQ CPU.

\subsection{Results}

\begin{center}
\begin{figure}
\includegraphics[scale=1]{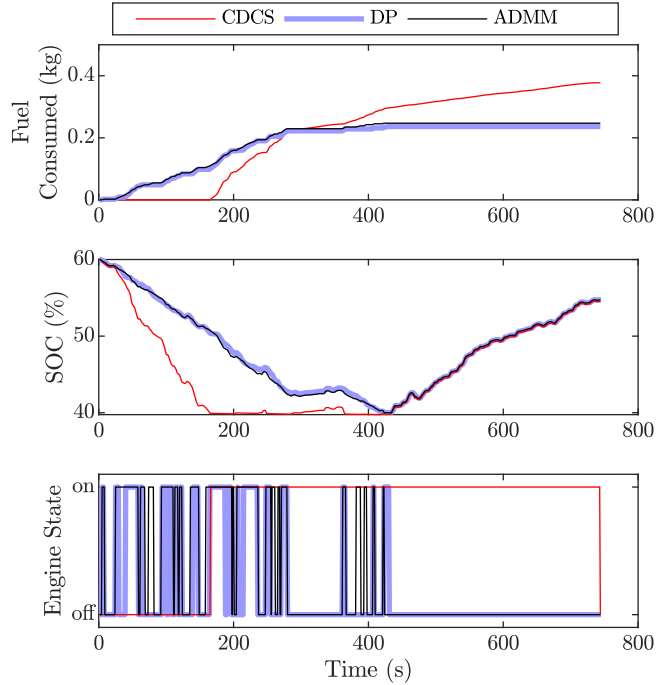}
\caption{Cumulative fuel consumption, SOC trajectory, and engine state for a single journey using CDCS, DP and ADMM.\label{fig_results2}}
\end{figure}
\end{center}

Figure \ref{fig_results2} shows the cumulative fuel consumption, state of charge trajectory, and engine switching control inputs for a single instance of the journey. It can be seen that although the state of charge (SOC) trajectories for DP and ADMM are qualitatively similar, the trajectory for CDCS has a large deviation for the first 430s, and this is reflected in the sub-optimal fuel consumption. After 430s the SOC trajectories are almost identical because the car is descending for the majority of the second half of the journey (see the gradient plots in Figure \ref{fig_velocity}), and is therefore in a regenerative mode for all three methods. During this period, the optimization-based methods consume no fuel as the engine has been switched off, but the CDCS strategy continues to consume fuel due to the rotation of the engine. It can also be seen from the engine controls that the ADMM algorithm largely obtains the periods where it is optimal to turn the engine on and off (e.g determines that the engine should be off while descending), but introduces additional switching decisions, particularly around 100s and 400s.

\begin{center}
\begin{figure}
\includegraphics[scale=1]{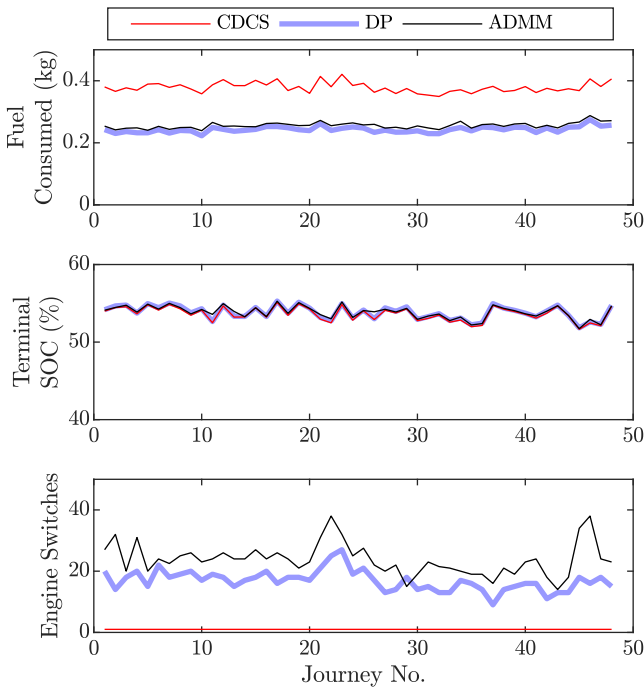}
\caption{Total fuel consumption, terminal SOC, and engine switching results for all journeys using CDCS, DP and ADMM. \label{fig_results1}}
\end{figure}
\end{center}

\begin{figure}
\begin{center}
\includegraphics[scale=1]{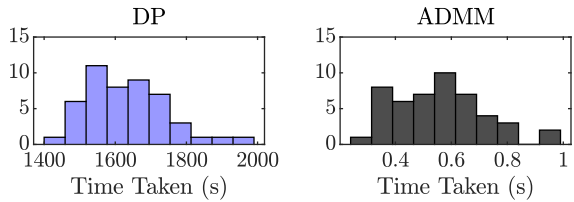}
\caption{Histograms of time taken for completion using DP and ADMM.\label{fig_results3}}
\end{center}
\end{figure}

The first plot in Figure \ref{fig_results1} shows the total fuel consumed across all journeys using CDCS, DP, and ADMM. To properly compare the fuel consumption using each method, the equivalent fuel consumption as a result of battery use is normally also considered within the total energy consumption, but in this case the second plot shows that the terminal SOC was within 1.7\% in all cases for all three methods, so this effect was ignored. It can be seen that DP reduces fuel consumption w.r.t CDCS by $\sim$40\% in all cases, and that whilst ADMM does not achieve the same level of optimality, on average it obtains 90.4\% of the fuel savings made using DP. It also shows that, on average, ADMM introduces 42.7\% more engine switches than DP. 
Figure \ref{fig_results3} shows histograms of the time taken to completion using DP and ADMM for all journeys. The mean time taken was 1,600s for DP and 0.56s for ADMM; a reduction of 99.97\%. It is worth noting that the algorithm is implemented here in vectorized Matlab code, so absolute performance could be increased further with a compiled software implementation where the $P_b$, $\sigma$, $\eta$ and $E$ variables are updated in parallel. Therefore, the results demonstrate that the ADMM algorithm is a promising candidate for a real-time shrinking/receding horizon MPC implementation.

\section{CONCLUSION}
In this paper we have demonstrated a locally convergent ADMM algorithm for the solution of the hybrid vehicle energy management problem, considering both power split and engine on/off decisions. Through simulations on real-world driving cycles, we have demonstrated that the ADMM algorithm can obtain 90.4\% of the fuel savings that are made when using a globally optimal DP algorithm, with a 3000-fold reduction in computational time. 

In this paper we have assumed that the future driver behaviour is known with complete precision; in reality it is subject to significant uncertainty. Our future work will investigate scenario-based approaches to address this limitation.

\section*{ACKNOWLEDGEMENT}
We would like to thank Professor Luigi del Re and Philipp Polterauer at Johannes Kepler University Linz, Austria, for generously sharing the driver data used in this publication and assisting with its processing.

\bibliographystyle{ieeetr}
\bibliography{references}

\end{document}